%% file: main.tex
\documentclass[letterpaper, 10 pt, conference]{ieeeconf}  

\IEEEoverridecommandlockouts                              
\usepackage{mathematix}
\usepackage{units}
\usepackage{flushend}
\usepackage{hyperref}
\hypersetup{
    bookmarks=true,         
    unicode=false,          
    pdftoolbar=true,        
    pdfmenubar=true,        
    pdffitwindow=false,     
    pdfstartview={FitH},    
    pdftitle={},    
    pdfauthor={},     
    pdfkeywords={},
    pdfsubject={Confidential},   
    pdfnewwindow=true,      
    colorlinks=true,        
    linkcolor=red,          
    citecolor=cyan,         
    filecolor=magenta,      
    urlcolor=cyan           
}

\usepackage{graphicx}      
\usepackage{amsmath,amsthm}
\usepackage{amssymb,gensymb}
\usepackage{xfrac}
\usepackage[latin1]{luainputenc}
\usepackage{eurosym}

\usepackage{acronym}

\input{acronyms.tex}

\newcommand{\costSingleShooting}{\phi}

\title{\LARGE \bf Aerial navigation in obstructed environments with
       embedded nonlinear model predictive control}

\author{Elias Small,
        Pantelis Sopasakis,
        Emil Fresk,
        Panagiotis Patrinos
        and George Nikolakopoulos
\thanks{
P. Sopasakis is with University of Cyprus,
Department of Electrical and Computer Engineering,
KIOS Research Center for Intelligent Systems and Networks,
1 Panepistimiou Avenue, 2109 Aglantzia, Nicosia, Cyprus.
{\tt\small sopasakis.pantelis@ucy.ac.cy}.
	}	%
\thanks{%
E. Small, E. Fresk and G. Nikolakopoulos are with the Robotics Team at Lule{\aa} Technical University,
Lule{\aa} SE-97187, Sweden
{\tt\small elias.small, emil.fresk, geonik@ltu.se}}%
\thanks{%
P. Patrinos is with the Department of Electrical Engineering (\textsc{Esat-Stadius}),
KU Leuven, Kasteelpark Arenberg 10, 3001, Leuven, Belgium.
{\tt\small panos.patrinos@esat.kuleuven.be}%
}
\thanks{This work was partially funded by the European Union's Horizon 2020
        Research and Innovation Programme under Grant Agreement No.730302 -- SIMS.
        P. Sopasakis was supported by European Union's Horizon 2020 research and
        innovation programme (KIOS CoE) under Grant No. 739551.
        P. Patrinos was supported by: FWO projects: No. G086318N; No. G086518N;
        Fonds de la Recherche Scientifique --- FNRS, the Fonds Wetenschappelijk Onderzoek
        --- Vlaanderen under EOS Project No. 30468160 (SeLMA) and
        Research Council KU Leuven C1 project No. C14/18/068.}
}
\begin{document}

\maketitle

\begin{abstract}
We propose a methodology for autonomous aerial navigation and obstacle avoidance
of \acp*{MAV} using \ac*{NMPC} and we demonstrate its
effectiveness with laboratory experiments. The proposed methodology can accommodate
obstacles of arbitrary, potentially non-convex, geometry.
The NMPC problem is solved using PANOC: a fast numerical optimization method which
is completely matrix-free, is not sensitive to ill conditioning, involves only
simple algebraic operations and is suitable for embedded NMPC. A \texttt{C89} implementation
of PANOC solves the NMPC problem at a rate of \(\unit[20]{Hz}\) on board a lab-scale
MAV. The MAV performs smooth maneuvers moving around an obstacle. For increased autonomy,
we propose a simple method to compensate for the reduction of thrust over time, which
comes from the depletion of the MAV's battery, by estimating the thrust constant.
\end{abstract}

\input{introduction.tex}

\input{uav_dynamics.tex}

\input{nmpc_formulation.tex}

\input{nmpc_solution.tex}

\input{experimental.tex}

\input{conclusions.tex}

\bibliographystyle{IEEEtran}
\bibliography{bibliography}

\end{document}

%% file: acronyms.tex
\acrodef{ode}       [ODE]           {Ordinary Differential Equation}
\acrodef{gas}       [GAS]           {Globally Asymptotically Stable}
\acrodef{ae}        [a.e.]          {almost everywhere}
\acrodef{wrt}       [w.r.t.]        {with respect to}
\acrodef{wlog}      [w.l.o.g.]      {without loss of generality}
\acrodef{AC}        [AC]            {absolutely continuous}
\acrodef{PMP}       [PMP]           {Pontryagin Maximum Principle}
\acrodef{lti}       [LTI]           {Linear Time-Invariant}
\acrodef{BLUE}      [BLUE]          {Best Linear Unbiased Estimator}
\acrodef{GP}        [GP]            {Gaussian Process}
\acrodef{MAP}       [MAP]           {Maximum A Posteriori}
\acrodef{MVUE}      [MVUE]          {Minimum Variance Unbiased Estimator}
\acrodef{MMSE}      [MMSE]          {Minimum Mean Square Error}
\acrodef{ML}        [ML]            {Maximum Likelihood}
\acrodef{LMMSE}     [LMMSE]         {Linear Minimum Mean Square Error}
\acrodef{LS}        [LS]            {Least Squares}
\acrodef{WSN}       [WSN]           {Wireless Sensor Network}
\acrodef{NARX}      [NARX]          {Nonlinear AutoRegressive eXogenus}
\acrodef{NCS}       [NCS]           {Networked Control System}
\acrodef{RKHS}      [RKHS]          {Reproducing Kernel Hilbert Space}
\acrodef{MSE}       [MSE]           {Mean-Square Error}
\acrodef{RMSE}      [RMSE]          {Root-Mean-Square Error}
\acrodef{RN}        [RN]            {Regularization Network}
\acrodef{CDF}       [CDF]           {Cumulative Distribution Function}
\acrodef{PDF}       [PDF]           {Probability Density Function}
\acrodef{UMP}       [UMP]           {Uniformly Most Powerful}
\acrodef{GLR}       [GLR]           {Generalized Likelihood Ratio}
\acrodef{UI}        [UI]            {Union-Intersection}
\acrodef{IU}        [IU]            {Intersection-Union}
\acrodef{SPRT}      [SPRT]          {Sequential Probability Ratio Test}
\acrodef{ISM}       [ISM]           {Industrial, Scientific and Medical}
\acrodef{DHT}       [DHT]           {Distributed Hash Table}
\acrodef{GPS}       [GPS]           {Global Positioning System}
\acrodef{UWB}       [UWB]           {Ultra Wide-Band}
\acrodef{RF}        [RF]            {Radio Frequency}
\acrodef{RFID}      [RFID]          {Radio Frequency Identification}
\acrodef{6LoWPAN}   [6LoWPAN]       {IPv6 over Low power Wireless Personal Area Networks}
\acrodef{AODV}      [AODV]          {Ad-hoc On-Demand Distance Vector}
\acrodef{ARW}       [ARW]           {Aerial Robotic Worker}
\acrodefplural{ARW} [ARWs] 	        {Aerial Robotic Workers}
\acrodef{MEMS}      [MEMS]          {Micro-Electro-Mechanical Systems}
\acrodef{IMU}       [IMU]           {inertial measurement unit}
\acrodef{UAV}       [UAV]           {Unmanned Aerial Vehicle}
\acrodefplural{UAV} [UAVs] 	        {Unmanned Aerial Vehicles}
\acrodef{MAV}       [MAV]           {micro aerial vehicle}
\acrodefplural{MAV} [MAVs] 	        {micro aerial vehicles}
\acrodef{VTOL}      [VTOL]          {Vertical Take-Off and Landing}
\acrodef{COR}       [CoR]           {Center of Rotation}
\acrodef{COG}       [CoG]           {Center of Gravity}
\acrodef{CAD}       [CAD]           {Computer Aided Design}
\acrodef{AR}        [AR]            {Auto Regressive}
\acrodef{MUSIC}     [MUSIC]         {Multiple Signal Classification}
\acrodef{KF}        [KF]            {Kalman Filter}
\acrodef{EKF}       [EKF]           {extended Kalman filter}
\acrodef{ESKF}      [ESKF]          {Error-State Kalman Filter}
\acrodef{SR-ESKF}   [SR-ESKF]       {Square-Root Error-State Kalman Filter}
\acrodef{UKF}       [UKF]           {Unscented Kalman Filter}
\acrodef{CKF}       [CKF]           {Cubature Kalman Filter}
\acrodef{MPC}       [MPC]           {model predictive control}
\acrodefplural{MPC} [MPC] 	        {model predictive controller}
\acrodef{DCM}       [DCM]           {Direction Cosine Matrix}
\acrodef{INS}       [INS]           {Inertial Navigation System}
\acrodef{SLAM}      [SLAM]          {Simultaneous Localization and Mapping}
\acrodef{SQP}       [SQP]           {sequential quadratic programming}
\acrodef{ROS}       [ROS]           {Robotic Operating System}
\acrodef{MSF}       [MSF]           {Multi Sensor Fusion}
\acrodef{PCA}       [PCA]           {Principal Component Analysis}
\acrodef{VIO}       [VIO]           {Visual Inertial Odometry}
\acrodef{TOF}       [TOF]           {Time of Flight}
\acrodef{RK}        [RK]            {Runge-Kutta}
\acrodef{NMPC}      [NMPC]          {non-linear model predictive control}
\acrodef{OS}        [OS]            {operating system}

\acrodef{MCU}       [MCU]           {Micro Controller Unit}
\acrodef{UART}      [UART]          {Universal Asynchronous Receiver Transmitter}
\acrodef{DSP}       [DSP]           {Digital Signal Processing}
\acrodef{CPU}       [CPU]           {Central Processing Unit}
\acrodef{FPU}       [FPU]           {Floating Point Unit}
\acrodef{DMA}       [DMA]           {Direct Memory Access}
\acrodef{SOC}       [SOC]           {System on a Chip}
\acrodef{PWM}       [PWM]           {Pulse Width Modulation}
\acrodef{PPM}       [PPM]           {Pulse Period Modulation}
\acrodef{ADC}       [ADC]           {Analog to Digital Converter}
\acrodef{I2C}       [I$^2$C]        {Inter-Integrated Circuit}
\acrodef{CAN}       [CAN]           {Controller Area Network}
\acrodef{USB}       [USB]           {Universal Serial Bus}
\acrodef{VCP}       [VCP]           {Virtual COM-Port}
\acrodef{ESD}       [ESD]           {Electro-Static Discharge}
\acrodef{GPS}       [GPS]           {Global Positioning System}
\acrodef{RSSI}      [RSSI]          {Received Signal Strength Indicator}
\acrodef{ESC}       [ESC]           {Electronic Speed Controller}
\acrodef{ISR}       [ISR]           {Interrupt Service Routine}
\acrodef{GPS}       [GPS]           {Global Positioning System}
\acrodef{BLDC}      [BLDC]          {Brushless Direct Current}
\acrodef{DC}        [DC]            {Direct Current}
\acrodef{PMSM}      [PMSM]          {Permanent Magnet Synchronous Motor}
\acrodef{RC}        [RC]            {Radio Control}
\acrodef{FCU}       [FCU]           {Flight Controller Unit}
\acrodefplural{FCU} [FCUs]          {Flight Controller Unit}

\acrodef{COBS}      [COBS]          {Consistent Overhead Byte Stuffing}
\acrodef{SLIP}      [SLIP]          {Serial Line Internet Protocol}
\acrodef{OTL}       [OTL]           {Over The Line}
\acrodef{PRE}       [PRE]           {Packet Pre-emption}

\acrodef{FROST}     [FROST]         {field robotics lab}
\acrodef{CARMA}     [CARMA]         {Compact AeRial MAnipulator}
\acrodef{LTU}		[LTU]			{Lule{\aa} University of Technology}
\acrodef{MFCQ}      [MFCQ]          {Mangasarian-Fromovitz constraint qualification}
\acrodef{LICQ}      [LICQ]          {linear independence constraint qualification}
\acrodef{IP}        [IP]            {interior point}
\acrodef{FBE}       [FBE]           {forward backward envelope}

%% file: introduction.tex
\section{Introduction}
\subsection{Background and motivation}
The need for safe aerial navigation and increased \ac{MAV}
autonomy nowadays poses
all the more relevant and pressing research questions, as drones make their
appearance in numerous application domains, such as the inspection of
critical or aging infrastructure \cite{Metni2007}, surveying of underground mines
\cite{Kanellakis+2019}, visual area coverage for search-and-rescue operations
\cite{Mansouri+2018}, precision agriculture \cite{Zhang2012} and many others.
In the majority of these applications, \acp{MAV} have to navigate in obstructed environments,
with static or moving obstacles of arbitrary geometry in known, or partially unknown surrounding environments.

Several methods have been proposed for navigation and collision
avoidance, such as potential field methods~\cite{minguezVPF2016,montiel2015bpf} and
graph search methods~\cite{ChengraphSearch2014}.
Alongside these methods, \ac{NMPC} is becoming popular for the navigation control
of various \acp{MAV} including fixed-wing aircrafts \cite{Stastny2017,Kamel2017}
and multi-rotor vehicles \cite{KAMEL20173463}.
\ac{NMPC} uses a nonlinear dynamical model of the system dynamics
to predict position and attitude trajectories from its current position
to a reference point, while avoiding all obstacles on its way and minimizing
a certain energy/cost function. In this way, a non-convex optimization problem needs
to be solved at every sampling time instant in a receding horizon fashion.
Another approach to obstacle avoidance is described in \cite{liu2017robust} where a high-level
path planner generates collision-free trajectories which are followed by an MPC controller.

In \cite{Son2017ModelPC}, \ac{SQP} is used to solve the \ac{NMPC} problem for
the navigation of a multi-rotor \ac{MAV} with a slung load, where the authors demonstrated
the effectiveness of \ac{NMPC}, however, provided neither evidence of the
solution quality or solver performance, nor an experimental verification.
\ac{NMPC} was used in \cite{Chao2012} for solving
obstacle and collision avoidance for several \acp{MAV} flying in formation, however again, only
simulations were done and the computation time was addressed.

Clearly, the presence of obstacle/collision avoidance constraints makes the
MPC problems particularly hard to solve.
\ac{SQP} is the method of choice in the literature \cite{Son2017ModelPC,Chao2012,Shim2006} that has as a main disadvantage the fact that it requires the solution of a
quadratic program (QP) at every iteration of the algorithm, which requires
inner iterations. SQP also requires computing and storing of the Jacobian matrices of the dynamics,
and sometimes the Hessians when the Hessian of the Lagrangian is used in the
QPs. Furthermore, the gradient descent method has been used to solve nonlinear MPC
problems for aerial navigation \cite{Shim2006}. This method, however, is sensitive to
bad conditioning and problems with long horizons tend to become
ill conditioned, while the convergence is expected to be slow.


\subsection{Contributions}
In this article we propose a control methodology for the autonomous
navigation of \acp{MAV} in obstructed environments. We allow for the obstacles
to have arbitrary non-convex shapes and, contrary to distance-based
methods~\cite{wang2014synthesis}, we do not require that the distance
function between the \acp{MAV} and each obstacle is available.

The \ac{NMPC} optimization problem is solved by using PANOC \cite{panocECC2018,SteThe+2017}
--- a recently proposed algorithm for non-convex optimization problems, which is
suitable for embedded \ac{NMPC}, as it requires only simple and cheap linear
operations (mainly inner products of vectors) and exhibits a fast convergence.
Unlike SQP, PANOC is matrix-free and only requires the computation of
Jacobian-vector products, which can be computed very efficiently by backward (adjoint)
automatic differentiation.
PANOC has been shown in \cite{panocECC2018,SteThe+2017,hermansPenalty2018}
to significantly outperform both \ac{SQP} and interior-point methods. To the authors best of knowledge, this is the first time that a fast \ac{NMPC} optimization problem is being demonstrated on an aerial platform, setting the base for future developments in the aerial robotics community. 

Our method for modeling has the strong merit of being independent of the mass of the
\ac{MAV}, whereas the norm in the community is to use mass and other detailed parameters of the
specific \ac{MAV} used, for example \cite{KAMEL20173463}, \cite{liu2017robust}, and \cite{Chao2012}.
This allows for our method to be used without tuning the specific mass or available thrust, improving
robustness, generalization and ease of use.

Evidence of the solution quality is provided by physical laboratory experiments, where a
\ac{MAV} is flown completely autonomously in a laboratory equipped with a VICON motion capture system.
The proposed method uses a full position and attitude model of the \acp{MAV}, which is able to run onboard, using 8-15\% CPU of a single core on an Intel
Atom Z8350.
As is shown in Section \ref{sec:experimental} the onboard controller is able to successfully
navigate the \ac{MAV} around an obstacle running at a sampling rate of \(\unit[20]{Hz}\)
and a prediction horizon of \(\unit[2]{s}\).

%% file: uav_dynamics.tex
\section{\ac{MAV} dynamics}
\subsection{\ac{MAV} kinematics}
The model of a quadrotor \ac{MAV}, defined by \cite{Kamel2017}, assumes that there exists
a low-level controller of  roll, pitch, yaw rate and thrust. This convention is 
common in \ac{MAV} flight controllers such as PixHawk, \cite{meier2012pixhawk}
and ROSFlight, \cite{rosflight}.  The high-level kinematics of the \ac{MAV} is given by
\begin{subequations}\label{eq:system-dynamics}
\begin{align}
        \dot{p}(t)
{}={}&
        v(t),
\\
        \dot{v}(t)
{}={}&
        R(\theta_r, \theta_p)
        \smallmat{0\\0\\T_d}
    {+}
        \smallmat{0\\0\\-g}
    -
        \smallmat{
            A_x & 0 & 0\\
            0 & A_y & 0\\
            0 & 0 & A_z}
        v(t),
\\
        \dot{\theta}_r(t)
{}={}&
        \nicefrac{1}{\tau_r}(K_r\theta_{r,d}(t)-\theta_r(t)),
\label{eq:system-dynamics-3}
\\
        \dot{\theta}_p(t)
{}={}&
        \nicefrac{1}{\tau_p}(K_p\theta_{p,d}(t)-\theta_p(t)),
\label{eq:system-dynamics-4}
\end{align}
\end{subequations}
where $p(t) = (p_x(t), p_y(t), p_z(t))\in \R^3$ and $v(t) \in \R^3$
are the position and velocity of the \ac{MAV} in the global frame of reference,
and $\theta_r \in \R$ and $\theta_p \in \R$ are the
roll and pitch angles, while $\theta_{r,d} \in \R$ and $\theta_{p,d} \in \R$ are
the reference angles sent to the low-level controller.
Furthermore, $T_d \in \Rplus$ is the $z$-axis thrust acceleration, while $A_x$, $A_y$, and
$A_z$ are the linear drag coefficients.
The lower layer --- the attitude control system --- is modeled by simple first-order
dynamics with time constant $\tau_r$ and $\tau_p$ and gains $K_r$ and $K_p$
for the roll and pitch respectively. Lastly, $R(\theta_r, \theta_p) \in \mathrm{SO}(3)$ describes the \ac{MAV}'s
attitude and is defined by the classical Euler angles in rotation matrix form as
\begin{align*}
 R(\theta_r, \theta_p) {}={} R_y(\theta_p)R_x(\theta_r),
\end{align*}
with
\begin{align*}
R_x(\theta_r) {}={}& \smallmat{
    1 & 0 & 0\\
    0 & \cos(\theta_r) & -\sin(\theta_r)\\
    0 & \sin(\theta_r) & \phantom{-}\cos(\theta_r)},
\\
R_y(\theta_p) {}={}& \smallmat{
    \phantom{-}\cos(\theta_p) & 0 & \sin(\theta_p)\\
    0 & 1 & 0\\
    -\sin(\theta_p) & 0 & \cos(\theta_p)}.
\end{align*}
Note that yaw is absent in this rotation matrix, as this model operates
in a yaw-compensated global frame, and the position control of the \ac{MAV} is
therefore independent of its yaw.
Moreover, it is important to note that we have chosen the acceleration, \(T_d\),
to be the manipulated variable of the system, rather than the corresponding
force, for the model to be mass-free.
This has the strong merit of making the controller robust to changes in the
mass of the \ac{MAV}, the available thrust from the motors, and the loss of
thrust over time due to the decline of battery voltage.

\subsection{Adaptive acceleration control}\label{secsub:thrust estimation}
\begin{figure}[t!]
    \includegraphics{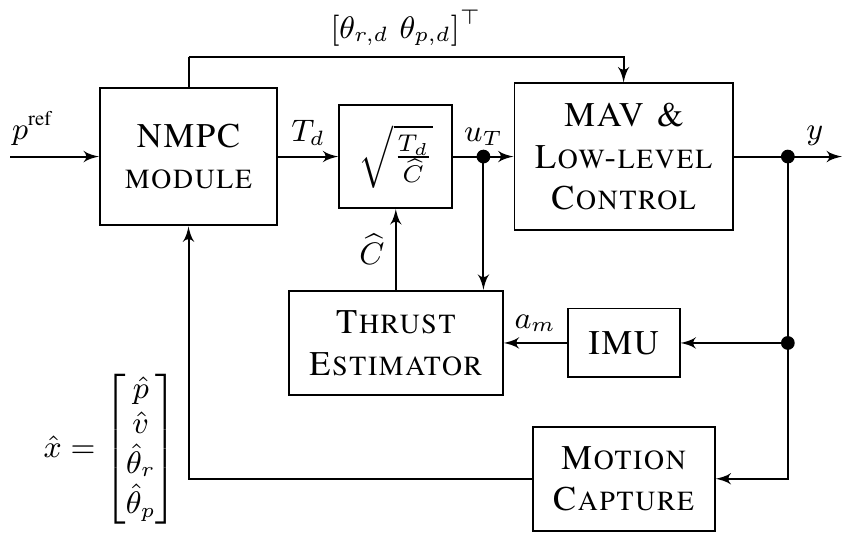}
    \vspace{-0.5cm}
    \caption{
        This diagram represents the complete \ac{MAV} system: \(p^\text{ref}\) is
	the reference position sent to the \ac{NMPC} controller, which
	calculates the desired angles, \(\theta_{r,d}\), \(\theta_{r,d}\), and thrust,
	\(T_d\). The thrust estimator uses the measured acceleration from the \ac{IMU},
	\(a_m\), to yield an estimate of the thrust constant, \(\widehat{C}\), which is
	used to obtain the thrust control signal \(u_T\). The complete system state,
    \(y\), is estimated by the motion capture system and measured by the \ac{IMU}
    which produces \(\hat{p},\hat{v},\hat{\theta}_r,\hat{\theta}_p\), which are
    sent to the \ac{NMPC} module, and the linear acceleration $a_m$ which is
    sent to the thrust estimator.}
    \label{fig:blockdiagram}
\end{figure}

In order for our design to be independent of the physical characteristics
that determine the available thrust acceleration, a simplified version of
\cite{fresk2018generalized} is used to continuously estimate the \ac{MAV}'s
maximum available thrust.
Following \cite{fresk2018generalized}, the force, \(F\), that is exercised
by the propellers, is given by
\begin{align}
    F = C_T u_T^2,
\end{align}
where $C_T$ is the \textit{thrust constant} and $u_T \in [0, ~1]$ is a unitless normalized thrust
control signal. The thrust constant is time dependent
based, for instance, on battery drain and how close the \ac{MAV} is to the ground, which
is why identifying a constant is not sufficient to track thrust references.
The issue is that there is no sensor in the system, which measures the generated force,
however, the \ac{IMU} can provide a measurement of the linear acceleration, albeit noisy.
Then, by dividing the thrust model by the mass of the \ac{MAV}, $m$, the model is now based
on acceleration:
\begin{subequations}
\begin{align}\label{eqn:thrust_acc}
    a = \frac{F}{m} = \frac{C_T}{m} u_T^2 = C u_T^2,
\end{align}
where $C \triangleq \nicefrac{C_T}{m}$ is the \textit{special thrust constant} of the vehicle.
Now the acceleration is measurable, together with
noise, and $u_T$ is what is sent to the low-level controllers, hence now it is
simply to choose the estimator of choice to estimate $C$.
Since $C$ is a slow moving parameter, we use the simple model
\begin{align}\label{eq:c_model}
    \dot{C} = \sigma_C^2 w,
\end{align}
where \(w\) is a zero-mean white noise%
\footnote{Equation \eqref{eq:c_model} is a stochastic differential equation which is
meant in the sense \(\dd C = \sigma_C^2 \dd B_t\), where \(B_t\) is the standard Brownian
motion.}.
\end{subequations}
Equations \eqref{eqn:thrust_acc} and \eqref{eq:c_model} define a nonlinear dynamical
system with state variable \(C\), input \(u_T\) and output \(y(C,u_T) = Cu_T^2\).
We estimate $C$ by means of an \ac{EKF}. \ac{EKF} is chosen because it is simple to 
tune, it allows to specify an initial estimate variance, and converges fast in the 
first few iterations.

In additional, we employ an outlier rejection scheme based on bounds of
the \emph{direct estimate} $\widetilde{C} \in [1g, ~10g]$, where
\begin{align}
    \widetilde{C} = \frac{a_m}{u_{T}^2},
\end{align}
which is calculated for each \ac{IMU} acceleration measurement $a_m$, which implies
that each acceleration measurement is inspected to enforce that no outliers are allowed
to update the filter.
These bounds result from the fact that a \ac{MAV} must be able
to generate at least \(1g\) of thrust to take off and it is assumed that
it cannot generate more than \(10g\) of thrust. 
The bounds on \(\widetilde{C}\) are inherited by the estimates \(\widehat{C}\)
yielding a simple constrained estimation scheme.

Once the thrust constant is estimated, an acceleration reference can be
converted to the thrust control signal $u_T$, by solving equation \eqref{eqn:thrust_acc} for $u_T$, resulting in
\begin{align}
    u_T = \sqrt{\frac{T_d}{\widehat{C}}},
\end{align}
A depiction of how the thrust constant estimation is tied to the overall scheme can be found in Fig.~\ref{fig:blockdiagram}.

\subsection{Overall system dynamics}\label{secsub:system dynamics}

The state of the controlled system is defined to be
\(x(t) = (p(t), v(t), \theta_r(t), \theta_p(t))\) and the manipulated input is
\(u(t) = (T_d(t), \theta_{r,d}(t), \theta_{p,d}(t))\). The system is observed using a VICON motion
capture system, which measures the full odometry of the system and provides the corresponding estimates
of the full state of the \ac{MAV} as $\hat{x} = (\hat{p}(t), \hat{v}(t), \hat{\theta}_r(t), \hat{\theta}_p(t))$.
Overall, the system dynamics can be concisely written as
\begin{equation}
    \dot{x}(t) = f(x(t), u(t)),
\end{equation}
where \(f\) is implicitly defined via \eqref{eq:system-dynamics}.

%% file: nmpc_formulation.tex
%
%
\section{Nonlinear MPC for obstacle avoidance}%
\label{sec:obstacle_avoidance}


\subsection{Navigation in obstructed environments}
We assume that a \ac{MAV} needs to navigate towards a reference position
\(p^{\mathrm{ref}} \in \R^3\), while avoiding a set of \(q(t)\in\N\) moving
obstacles \(\{O_{j}(t)\}_{j\in\N_{[1,q(t)]}}\).

We select \(n_{\mathrm{mav}}\) \emph{corner points} on the \ac{MAV} and
position a ball with radius \(r_{\mathrm{ball}}\) centered at each such point
so that the whole vehicle is contained in the union of these balls.
We assume that the coordinates of the corner points in the global frame of
reference are given by \(c_{\iota}(p(t))\), for \(\iota\in\N_{[1,n_{\mathrm{mav}}]}\).

\begin{figure}
    \centering
    \includegraphics[width=0.8\linewidth]{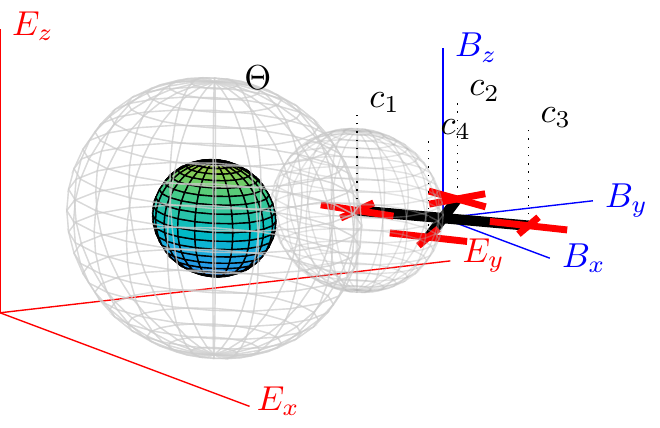}
    \caption{A quadrotor and an spherical obstacle \(O(t)\) (colored solid ball)
    and its enlargement \(\Theta(t)\). We have selected four corner points, 
    \(c_1, c_2,c_3,c_4\) on the \ac{MAV}. The red lines indicate the earth-fixed 
    frame of reference, \((E_x, E_y, E_z)\), and the blue ones the body-fixed frame, 
    \((B_x, B_y, B_z)\).}
    \label{fig:enlarged-obstacle}
\end{figure}

In order for the \ac{MAV} to not collide with the obstacles, we shall require
that
\begin{equation}
    c_{\iota}(p(t))
{}\notin{}
    \Theta_{j}(t)
{}\dfn{}
    O_{j}(t) + \mathcal{B}_{r}(r_{\mathrm{ball}}),
\end{equation}
for all \(j\in\N_{[1, q(t)]}\), \(\iota\in\N_{[1,n_{\mathrm{mav}}]}\), where
\(\mathcal{B}_{r}(r_{\mathrm{ball}})\) is a ball centered at the origin
with radius \(r_{\mathrm{ball}}\).
The set \(\Theta_{j}(t)\) is an enlarged version of the original
obstacle \(O_{j}(t)\).
The concept is illustrated in Fig.~\ref{fig:enlarged-obstacle}.

We introduce the stage cost function
\(
	\ell
{}:{}
	\R^{n_x}
    {}\times{}
	\R^{n_u}
    {}\times{}
	\R_+
    {}\to{}
	\R_+
\)
and the terminal cost function
\(
	\ell_{f}
{}:{}
	\R^{n_x}
    {}\times{}
	\R_+
    {}\to{}
	\R_+
\) which penalize the deviation of the system state from a reference state.
Typical choices are
\begin{subequations}
\begin{align}
        \ell(x, u, t)
{}={}&
        \|x-x^{\mathrm{ref}}(t)\|_{Q}^{2}
    {}+{}
        \|u-u^{\mathrm{ref}}(t)\|_{R}^{2},
\\
        \ell_f(x, t)
{}={}&
        \|x-x^{\mathrm{ref}}(t)\|_{Q_f}^2,
\end{align}
\end{subequations}
where \(Q\in\R^{n_x\times n_x}\), \(R\in\R^{n_u\times n_u}\) and 
\(Q_f\in\R^{n_x\times n_x}\) are positive semi-definite matrices and
\(x^{\mathrm{ref}}\) is the reference state which has the form
\(x^{\mathrm{ref}} = [p^{\mathrm{ref}} ~ 0_{1\times n_x - 3}]^{\top}\).

The nonlinear model predictive control problem for navigation in an obstructed
environment consists in solving the following problem
\begin{subequations}\label{eq:nmpc-1}
\begin{align}
 \minimize\ &J {\,}{=}{\,} \ell_{f}(\bar{x}(T), T) {\,}{+}\hspace{-0.1em}
 \int_{0}^{T}\hspace*{-0.7em}\ell(\bar{x}(\tau), \bar{u}(\tau), \tau)\dd \tau \label{eq:nmpc-1-1}\\
 \stt\ &\bar{x}(0) = \hat{x},\label{eq:nmpc-1-2}\\
 &\dot{\bar{x}}(t) = f(\bar{x}(t), \bar{u}(t)),\label{eq:nmpc-1-3}\\
 &\bar{u}(t) \in U(t),\label{eq:nmpc-1-4}\\
 &c_{\iota}(\bar{p}(t))\notin \Theta_{j}(t),\, j\in\N_{[1,q(t)]},\label{eq:nmpc-1-5}\\
 &\qquad \iota \in \N_{[1, n_{\mathrm{uav}}]},\ t\in[0,T]\notag
\end{align}
\end{subequations}
where \(\bar{u}(t) = (\bar{T}_d(t), \bar{\theta}_{r,d}(t), \bar{\theta}_{p,d}(t))\) and 
\(\bar{x}(t) = (\bar{p}(t), \bar{v}(t),\) \( \bar{\theta}_{r}(t), \bar{\theta}_{p}(t))\), 
for \(t\in[0,T]\) are the predicted input and state signals.

In this formulation we have assumed that the future trajectories of all
obstacles are exactly known and independent of the trajectory of the
controlled vehicle. If this is not the case, we have to formulate
appropriate robust or stochastic variants of the above obstacle
avoidance problem.

The control action is exercised to the system via a zero-order hold
element, that is, \(\bar{u}(t) = \bar{u}_{k}\) for \(t\in[kT_s, (k+1)T_s)\),
where \(T_s\) is the sampling period. We assume that \(T=NT_s\) for some
\(N\in\N\). Then, the cost function in \eqref{eq:nmpc-1-1} can be written as
\begin{equation}
    J = \ell_{f}(\bar{x}(T), T) + \sum_{k=0}^{N-1}\int_{kT_s}^{(k+1)T_s}
    \hspace*{-1.8em}\ell(\bar{x}(\tau), \bar{u}_k, \tau)\dd \tau.
\end{equation}
Since it is not possible to derive analytical solutions of the
nonlinear dynamical system \eqref{eq:nmpc-1-3}, the system trajectories
as well as the cost function \(J\) along these trajectories has to
be evaluated by discretizing the system dynamics and integrals.
By doing so, the system state trajectoriy \(\bar{x}(t)\) is
evaluated at points \(\bar{x}_k = \bar{x}(kT_s)\) as follows
\begin{equation*}
    \bar{x}_{k+1} \approx f_k(\bar{x}_k, \bar{u}_k),
\end{equation*}
and
\begin{equation*}
    J \approx \ell_N(\bar{x}_N) + \sum_{k=0}^{N-1}\ell_k(\bar{x}_k, \bar{u}_k).
\end{equation*}
Any explicit integration method such as the fourth-order Runge-Kutta or Forward
Euler lead to high quality approximations of \ac{MAV} trajectories.
This way, the original continuous-time optimal control problem is
approximated by a discrete-time one which is solved at every time
instant in a receding horizon fashion.

\subsection{Penalty functions for obstacles of general shape}
Each obstacle is described by a set of \(m_j(t)\) nonlinear constraints
of the form
\begin{equation}
    \Theta_{j}(t) =
        \{
            p\in\R^{3} \mid h_{j}^{i}(p, t) {}>{} 0,
            i\in\N_{[1, m_j(t)]}
        \},
\end{equation}
where functions \(h_{j}^{i}{}:{}\R^3\times\R_+\to\R_+\) are
\(C^{1,1}\) functions.
This approach allows one to describe obstacles of very general convex or
nonconvex shape.
For example, by choosing functions \(h_j^i\) to be affine in \(p\), we
can model any polytopic object. Functions of the form
\(h_j(p,t) = 1 - (p-p_0(t))^{\top}M(t)(p-p_0(t))\) can be used to model
ellipsoidal objects or elliptic cylindrical ones. Polynomial, trigonometric
and other functions can be used to model more complex geometric shapes.

\begin{figure}
    \centering
    \includegraphics[width=0.48\linewidth]{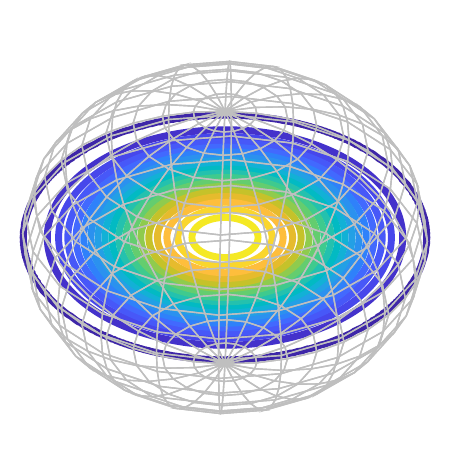}
    \includegraphics[width=0.40\linewidth]{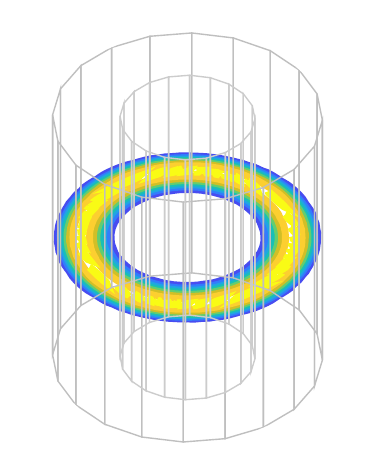}
    \caption{Level sets of slices of the function \(\psi_{\Theta(t)}\)
    on the plane \(\{(p_x,p_y,p_z)\in\R^3 {}\mid{} p_z=0\}\)
    for (Left) a ball-shaped obstacle and (Right) a non-convex obstacle.
    The obstacles are circumscribed by light gray mesh lines.}
    \label{fig:penalty_function}
\end{figure}

For simplicity, in this section we focus in the case where there is
one obstacle, that is \(q(t)=1\), which we denote by
\(
        \Theta(t)
{}={}
        \{
                p\in\R^3
            {}\mid{}
                h^i(p, t) > 0,
                i\in\N_{[1,m]}
        \}
\).
The constraint \(c_\iota(p(t))\notin O(t)\) is satisfied if and only if
\(h^{i_0}(c_\iota(p(t)), t) \leq 0\) for some \(i_0\in\N_{[1,m]}\), or
equivalently, if
\begin{equation}
    \psi_{\Theta(t)}(c_\iota(p(t))) {}={} 0,
\end{equation}
for all \(\iota \in \N_{[1, n_{\mathrm{uav}}]}\), where \(\psi_{\Theta(t)}:\R^3\to\R_+\)
is the function defined as
\begin{equation}
        \psi_{\Theta(t)}(p)
{}\dfn{}
        \tfrac{1}{2}\prod_{i=1}^{m}[h^{i}(p, t)]_+^2.
\end{equation}
Such a function is illustrated in Fig. \ref{fig:penalty_function}. Function
\(\psi_{\Theta(t)}(p)\) takes the value \(0\) outside the enlarged obstacle
\(\Theta(t)\) and increases in the interior of it as we move away from its boundary.

Function \(\psi_{\Theta(t)}\) is differentiable with gradient
\[
	\nabla \psi_{\Theta(t)} (p)
{}={}
	1_{\Theta(t)}(p)
	\sum_{i=1}^{m} h^i(p, t)
	  \prod_{j\neq i}
	    (h^{j}(p, t))^2\nabla_{p} h^i(p, t),
\]
where \(1_{\Theta(t)}\) is the characteristic function of \(\Theta(t)\) with
\(1_{\Theta(t)}(p) {}={} 1\) if \(p\in\Theta(t)\) and \(1_{\Theta(t)}(p) {}={} 0\)
otherwise.

Functions \(\psi_{\Theta(t)}\) can be used to impose the obstacle avoidance
requirements as soft constraints. To this end, we eliminate the non-convex
constraints \(c_{\iota}(\bar{p}(t)) \notin \Theta_j(t)\) and introduce the
modified stage and terminal cost functions
\begin{subequations}\label{eq:nmpc_cost_functions}
\begin{align}
    \tilde{\ell}(\bar{x}, \bar{u}, t)
{}={}&
    \ell(\bar{x}, \bar{u}, t)
{}+{}
    \sum_{\iota, j} \lambda_{j, \iota}
    \psi_{\Theta_{j}(t)}(c_{\iota}(\bar{p}(t))),
\\
    \tilde{\ell}_f(\bar{x}, \bar{u}, t)
{}={}&
    \ell_f(\bar{x}, \bar{u}, t)
{}+{}
    \sum_{\iota, j} \lambda^{f}_{j, \iota}
      \psi_{\Theta_{j}(t)}(c_{\iota}(\bar{p}(t))),
\end{align}
\end{subequations}
where \(\lambda_{j, \iota}\) and \(\lambda^{f}_{j, \iota}\) are positive weight
coefficients. The overall \ac{MPC} problem becomes
\begin{subequations}\label{eq:nmpc_discrete_time}
 \begin{align}
  \minimize_{\bar{u}_0,\ldots, \bar{u}_{N-1}}{} \
      &
      \tilde{\ell}_{N}(\bar{x}_N)
  {}+{}
      \sum_{k=0}^{N-1} \tilde{\ell}_{k}(\bar{x}_k, \bar{u}_k)
\\
 \stt\
      & \bar{x}_0 = x
\\
      & \bar{x}_{k+1} = f_k(\bar{x}_k, \bar{u}_k),\ k\in\N_{[0,N-1]}
      \label{eq:nmpc_discrete_time:dynamics}
\\
      & \bar{u}_k \in U_k,\ k\in\N_{[0,N-1]}
 \end{align}
where \(U_k = U(kT_s)\).
\end{subequations}
The optimization is carried out over finite-dimensional vectors
\(\bar{u} = (\bar{u}_0, \ldots, \bar{u}_{N-1})\in\R^{n}\) with
\(n {}={} n_u(N-1)\).

\subsection{Single-shooting problem formulation}
We shall cast optimization problem \eqref{eq:nmpc_discrete_time} in the
following compact and simple form
\begin{equation}\label{eq:nmpc_compact}
 \minimize_{\bar{u}\in U}{} \costSingleShooting(\bar{u};\, \hat{x}, p^{\mathrm{ref}}({}\cdot{})),
\end{equation}
where
\(
	U
{}={}
	U_0
    {}\times{}
	U_1
    {}\times{}
	\ldots
    {}\times{}
	U_{N-1}
\) and
\(
	\costSingleShooting
{}:{}
	\R^{n}
{}\to{}
	\R_{+}
\)
is a \(C^{1,1}\) function. To this end, we need to eliminate the state
sequence in \eqref{eq:nmpc_discrete_time:dynamics}. Let us introduce a
sequence of functions \(F_{k} {}:{} \R^{n} {}\to{} \R^{n_x}\) for
\(k\in\N_{[0, N]}\) defined recursively by
\begin{subequations}
 \begin{align}
  F_0(\bar{u}) {}={}& \hat{x},
  \\
  F_{k+1}(\bar{u}) {}={}& f_k(F_{k}(\bar{u}), \bar{u}_k).
 \end{align}
\end{subequations}
Then, problem \eqref{eq:nmpc_discrete_time} is written as in
\eqref{eq:nmpc_compact} with
\begin{equation}
	\costSingleShooting(\bar{u})
{}={}
	\tilde{\ell}_{N}(F_N(\bar{u}))
{}+{}
	\sum_{k=0}^{N-1} \tilde{\ell}_{k}(F_k(\bar{u}), \bar{u}_k).
\end{equation}
This is known as the \textit{single shooting} formulation \cite{panocECC2018}.

%% file: nmpc_solution.tex
\section{Fast online nonlinear \ac{MPC} using PANOC}

Problem \eqref{eq:nmpc_compact} is in a form that can be solved by PANOC 
\cite{panocECC2018}. In particular, the gradient of \(\costSingleShooting\) can be computed using 
automatic differentiation \cite{dunn1989efficient} which is implemented by software 
such as CasADi \cite{Andersson2018}. PANOC finds a \(\bar{u}^\star\in\R^{n}\) which solves the 
optimality conditions 
\begin{equation}\label{eq:zero_fpr}
 R_{\gamma}(\bar{u}^\star) {}={} 0,
\end{equation}
where \(R_{\gamma}:\R^{n} {}\to{} \R^{n}\) is the \textit{fixed-point 
residual} operator with parameter \(\gamma{}>{}0\) defined as 
\begin{equation}
 R_{\gamma}(\bar{u}) {}={} \bar{u} - T_{\gamma}(\bar{u}),
\end{equation}
where \(T_\gamma:\R^{n} {}\to{} \R^{n}\) is the \textit{projected gradient 
operator} given by
\begin{equation}
 T_{\gamma}(\bar{u}) {}={} \proj_{U}(\bar{u} - \gamma \nabla \costSingleShooting(\bar{u})).
\end{equation}

PANOC combines \textit{safe} projected-gradient updates \(\bar{u}^{\nu+\nicefrac{1}{2}}\)
with \textit{fast} Newton-type directions \(d^{\nu}\) computed by L-BFGS while it uses
the \ac{FBE} function \(\varphi_{\gamma}\) as a merit function for globalization given by
\begin{equation}
    \varphi_{\gamma}(\bar{u}) 
{}={}
    \costSingleShooting(\bar{u}) 
{}-{} 
    \nicefrac{\gamma}{2} \|\nabla \costSingleShooting(\bar{u})\|^2 
{}+{} 
    \dist_{U}^2(\bar{u} {}-{} \gamma \nabla \costSingleShooting(\bar{u})).
\end{equation}
The forward-backward envelope is an exact, continuous and real-valued merit
function which shares the same (local/strong) minima with \eqref{eq:nmpc_compact}.
That said,  Problem \eqref{eq:nmpc_compact} is reduced to the unconstrained minimization 
of \(\varphi_\gamma\).

PANOC is shown in Algorithm \ref{alg:panoc}.
L-BFGS uses a buffer of length \(\mu\) of vectors 
\(
	s^{\nu} 
{}={} 
	\bar{u}^{\nu+1} {}-{} \bar{u}^{\nu}
\) and 
\(
	y^\nu 
{}={} 
	R_{\gamma}(\bar{u}^{\nu+1}) 
{}-{} 
	R_\gamma(\bar{u}^{\nu})
\) to compute the update directions \(d^{\nu}\)~\cite{panocECC2018},{} \cite[Sec.~7.2]{nocedal2006numerical}. 
The computation of \(d^{\nu}\) requires only inner products which amount to a maximum 
of \(4\mu n\) scalar multiplications.
In particular, following \cite{Li2001}, the L-BFGS buffer is updated only if 
 \(
 	s^{\nu\top}y^{\nu}/\|s^{\nu}\|^2 
     {}>{} 
 	\epsilon_{d} \|R_{\gamma}(\bar{u})\|
 \).

Overall, PANOC uses exactly the same oracle as the projected gradient method, that is 
it only requires the invocation of \(\proj_{U}\), \(\costSingleShooting\) and 
\(\nabla \costSingleShooting\).
Lastly, owing to the \ac{FBE}-based line search, PANOC converges globally, that is, 
from any initial guess, \(\bar{u}^0\).

\begin{algorithm}
\caption{PANOC algorithm for nonlinear \ac{MPC}}\label{alg:panoc}
 \input{panoc.tex}
\end{algorithm}


%% file: panoc.tex
\newcommand{\ubarnu}{\bar{u}^{\nu}}
\newcommand{\ubarnuplus}{\bar{u}^{\nu+1}}
\newcommand{\ubarhalf}{\bar{u}^{\nu+\nicefrac{1}{2}}}

\begin{algorithmic}
 \REQUIRE 	 	
 	Initial guess \(\bar{u}^0\in\R^{n}\),\
 	Current state \(x\in\R^{n_x}\),
 	Estimate \(L>0\) of the Lipschitz constant of \(\nabla \costSingleShooting\),\
 	L-BFGS memory length \(\mu\),
 	Tolerance \(\epsilon>0\),
 	Maximum number of iterations \(\nu_{\mathrm{max}}\)%
 \ENSURE Approximate solution \(\bar{u}^{\star}\) 	
\STATE Choose \(\gamma\in(0,\nicefrac{1}{L})\), \(\sigma\in(0,\tfrac{\gamma}{2}(1-\gamma L))\)
\FOR{\(\nu=0,1,\ldots,\nu_{\mathrm{max}}\)}
 	\STATE
 		Compute \(\nabla\costSingleShooting(\bar{u}^\nu)\) using
 		automatic differentiation
	\STATE\label{state:zerofpr2:FB}
		\(
			\ubarhalf
		{}\gets{}
			T_{\gamma}(\bar{u}^\nu)
		\)
	\STATE	\label{state:zerofpr2:r}
		\(
			r^{\nu}
		    {}\leftarrow{}
			\gamma^{-1}(\bar{u}^\nu{}-{}\ubarhalf)
		\)
	\STATE \algorithmicif{} \(\|r^{\nu}\| {}<{} \epsilon\), exit
	\WHILE {\(
		      \costSingleShooting(\ubarhalf)
		 {>}
		      \costSingleShooting(\ubarnu) 
		 {-}
		      \gamma \nabla \costSingleShooting(\ubarnu)^{\top}r^{\nu}
		 {+} \tfrac{L}{2}\|\gamma r^{\nu}\|^2
		\)}
	  \STATE Empty the L-BFGS buffers
	  \STATE \(L {}\gets{} 2L\), 
		 \(\sigma {}\gets{} \nicefrac{\sigma}{2}\), 
		 \(\gamma {}\gets{} \nicefrac{\gamma}{2}\)
	  \STATE \(\ubarhalf {}\gets{} T_{\gamma}(\ubarnu)\)
	\ENDWHILE	
 	\STATE \label{state:zerofpr2:d} \(d^\nu\leftarrow-H_\nu r^\nu\) using L-BFGS%
	\STATE \label{state:zerofpr2:upd} 
		\(
			\ubarnuplus
		    {}\gets{}
			\ubarnu 
		    {}-{} 
			(1-\tau\!_\nu) \gamma r^\nu + \tau\!_\nu d^\nu
		\),
		where \(\tau\!_\nu\) is the largest number in \(\{\nicefrac{1}{2^i}{}:{}i\in\N\}\) such that
		\(
			\varphi_\gamma(\ubarnuplus)
		{}\leq{}
			\varphi_\gamma(\ubarnu){}-{}\sigma\|r^\nu\|^2
		\)
\ENDFOR{}%
\end{algorithmic}

%% file: experimental.tex
\section{Experimental validation}\label{sec:experimental}
\begin{figure}[thb!]
\centering
      \includegraphics[width=0.8\linewidth]{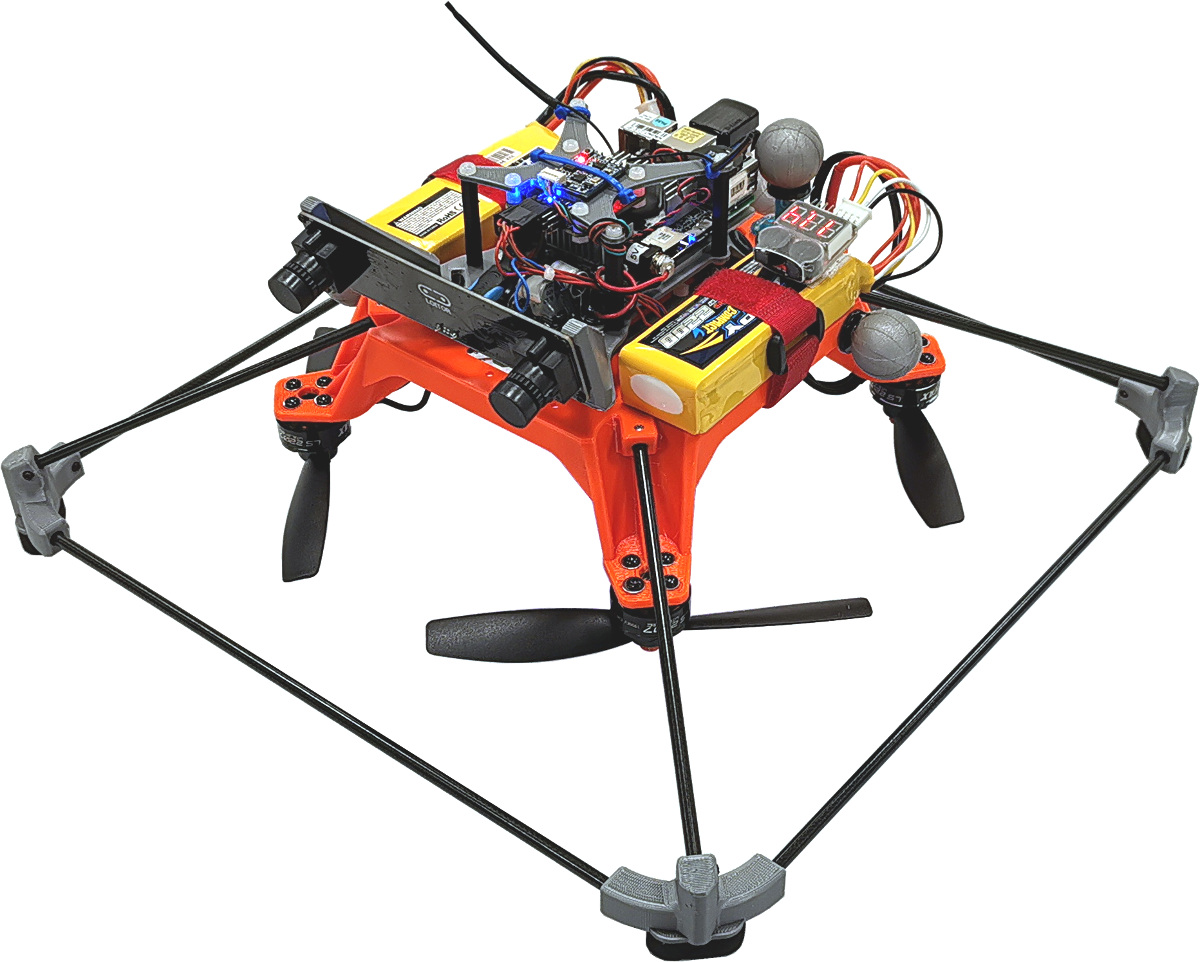}
      \caption{The inverted quadrotor used in the experiments, which is
               specifically designed to have a small $x/y$ footprint of \(\unit[34]{cm}\) by \(\unit[34]{cm}\), a height of \(\unit[12]{cm}\), and weight of \(\unit[1.02]{kg}\), to be
               suitable for indoor flight.}%
      \label{fig:invquad}
\end{figure}

For the experimental validation of the proposed control scheme, an inverted
quadrotor using the ROSFlight \cite{rosflight} low-level controller was used
for all trials, shown in Fig.~\ref{fig:invquad}. The onboard computer used is
an Aaeon Up Board with an Intel Atom x5-z8350 processor with four
\(\unit[1.44]{GHz}\) cores and \(\unit[2]{GB}\) of RAM.
The board runs Ubuntu Server 16.04. The field robotics lab at Lule\aa University of technology is equipped with a Vicon motion capture system featuring 20 infrared cameras that track the odometry of the \ac{MAV};
this data is used by the \ac{NMPC} controller for navigation.

The NMPC module runs simple \texttt{C89} code which was generated by
\texttt{nmpc-codegen} --- an LGPLv3.0-licensed open-source code generation
toolkit which is available at \url{github.com/kul-forbes/nmpc-codegen}.

An upright cylindrical obstacle, \(O\), is placed so that its vertical symmetry
axis runs through the origin \((0.0, 0.0, 0.0)\) of the global coordinate frame
in the flying arena at \ac{FROST}. The cylinder, \(O\),  has a radius of
\(r_{\mathrm{cyl}} = \unit[0.45]{m}\) and height \(z_{\mathrm{cyl}} {}={} \unit[2]{m}\). 
The obstacle is described by the functions
\(
	h^1(p, t)
{}={}
	r^2-p_x^2-p_y^2
\),
\(
	h^2(p, t)
{}={}
	p_z
\)
and
\(
	h^3(p, t)
{}={}
	z_{\mathrm{cyl}} - p_z
\).
A single corner point is used which is positioned at the center of the \ac{MAV};
the enclosing ball as in Fig. \ref{fig:enlarged-obstacle} has a radius of 
\(r_{\mathrm{ball}} {}={} \unit[0.24]{m}\).
In order to account for possible small constraint violations due to the fact
that obstacle avoidance constraints are modeled via penalty functions,
we consider an additional enlargement of \(\unit[0.06]{m}\).
As a result, the enlarged cylinder, \(\Theta(t)\), has a radius of
\(\unit[0.75]{m}\) and height \(\unit[2.3]{m}\).
The weights of the obstacle constraints, $\lambda_{j, \iota}$ and
\(\lambda^{f}_{j, \iota}\), in Equation~\eqref{eq:nmpc_cost_functions} were all 
set to $10000$, and the continuous-time system was integrate with the forward
Euler method.

\begin{figure}[!ht]
\centering
      \includegraphics{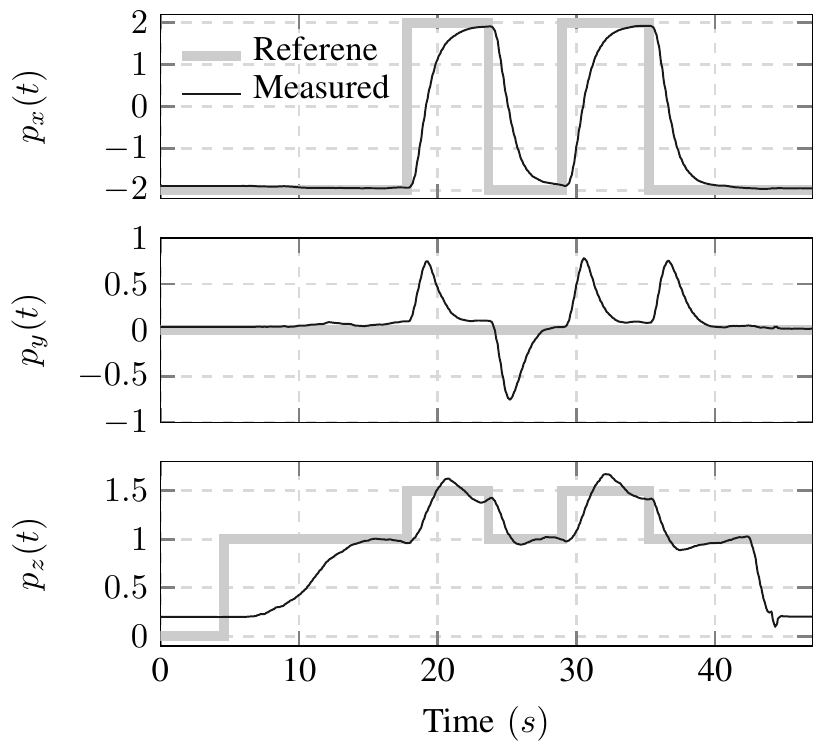}
      \caption{Smooth navigation of the \ac{MAV} in space: the position of the
               vehicle versus time. Positions are in \(\unit{m}\) and time is in \(\unit{s}\).}
      \label{fig:xyz-path}
\end{figure}

The flight test performed for avoiding the obstacle consisted
of alternating between two position references on opposite sides of the obstacle.
The two position references given alternately were, in meters,
\((-2.0, 0.0, 1.0)\) and \((2.0, 0.0, 1.5)\). These references were sent when
the \ac{MAV} was close to its previous reference position. The exact time the references are
changed can be seen in Fig.~\ref{fig:xyz-path}.

\ac{NMPC} runs at \(\unit[20]{Hz}\) with a prediction and control horizon of 40 steps, meaning the solver predicts
the states of the system \(\unit[2]{s}\) into the future.
The solver occupied between $8\%$ and $15\%$ of CPU on an Intel Atom Z8350 --- 
an indication of the solver's computational efficiency. 

\begin{figure}[!ht]
\centering
      \includegraphics{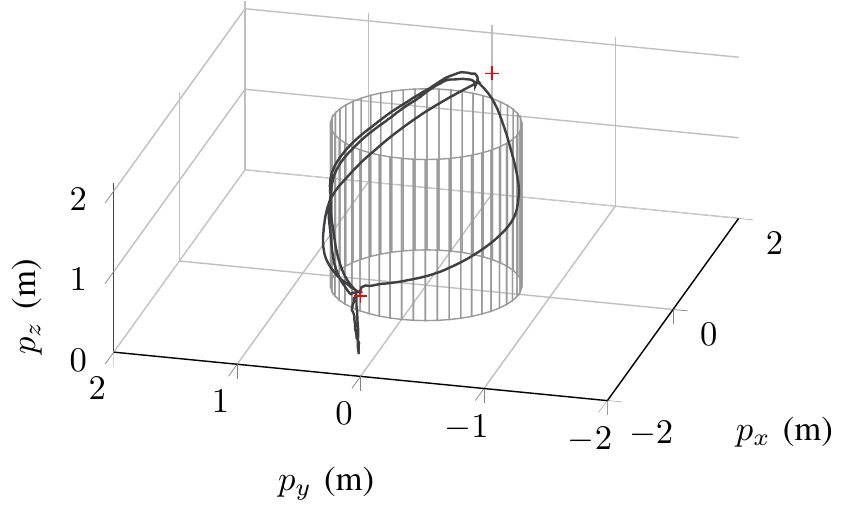}
      \caption{Path of the UAV autonomously taking off, traveling between the two
              reference positions \(p^{\mathrm{ref}, 1} {}={} (-2, 0, 1)\) and
              \(p^{\mathrm{ref}, 2} {}={} (2, 0, 1.5)\), and landing. An upright cylindrical
              (enlarged) obstacle of radius \(\unit[75]{cm}\) is positioned so that its axis
              runs through \((0,0,0)\).}
      \label{fig:path}
\end{figure}

Fig.~\ref{fig:path} shows the actual path flown by the \ac{MAV} during the test where the positioning
data is taken from the motion tracking system and has sub-millimeter accuracy. 
The path is also shown in Fig.~\ref{fig:xyz-path} where we plot the \ac{MAV}'s position versus time.
The \ac{MAV} does not have time to settle at the
reference altitude as a new reference is sent to the controller before the position completely converges.

\begin{figure}[!ht]
      \includegraphics{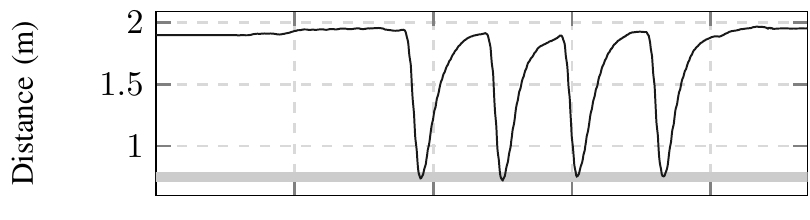}
      \caption{Distance of the center of the MAV from the center of the
               enlarged cylinder, \(\Theta\). The solid grey line indicates the
               border of the cylinder at a radius of \(\unit[75]{cm}\).}

      \label{fig:distancetocenter}
\end{figure}

As the \ac{MAV} passes the obstacle it violates the obstacle constraint, as shown in
Fig.~\ref{fig:distancetocenter}, which is expected from the penalty formulation.
The maximum violation is \(\unit[2.86]{cm}\), which is below the extra enlargement
of \(\unit[6]{cm}\) of the obstacle.

\begin{figure}[!ht]
      \includegraphics{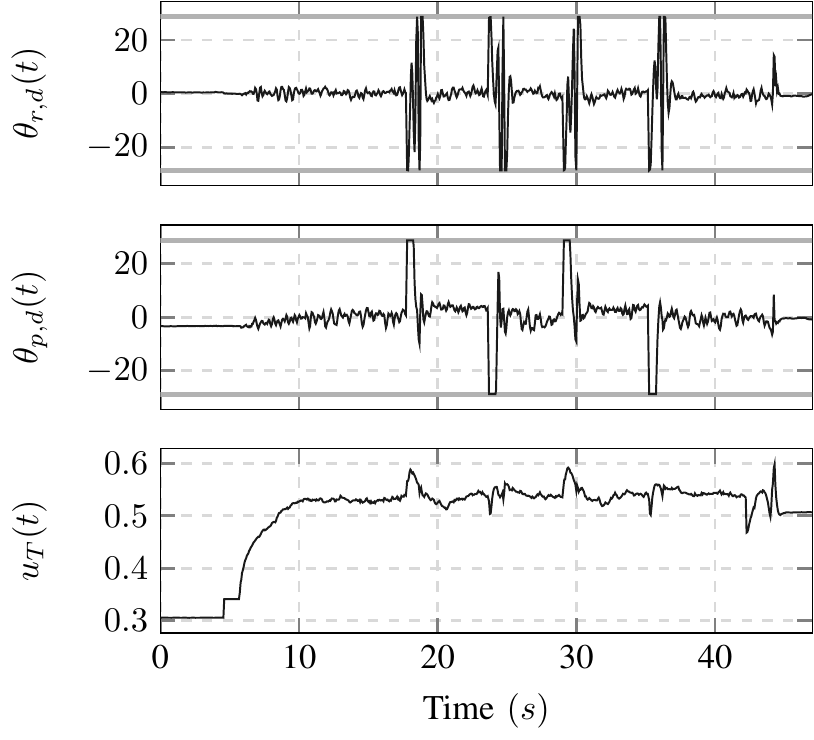}
      \caption{Control signals sent from the solver to the low-level controller 
               during the experiment. The angles are in degrees for ease of reading.}
      \label{fig:controlsignals}
\end{figure}

Fig.~\ref{fig:controlsignals} shows the control signals (roll, pitch, and normalized thrust
references) commanded by the \ac{NMPC}. The roll and pitch angles are bound between \(\unit[-0.5]{rad}\)
and \(\unit[0.5]{rad}\); these bounds are active as shown in Fig.~\ref{fig:controlsignals}.
This further motivates the use of \ac{NMPC}, allowing for bounds to be
directly included in the problem formulation. 

The control signals could be made less aggressive by penalizing the rate of change of the input
in \eqref{eq:nmpc_discrete_time}, that is, by adding a penalty of the form 
\(
	\ell_{\Delta} 
{}={}
	\|\bar{u}_k - \bar{u}_{k-1}\|^2_{R_{\Delta}}
\)
for a symmetric positive semidefinite matrix \(R_{\Delta} \in \R^{n_u\times n_u}\).
Nevertheless, the maneuvering of the \ac{MAV} is smooth as shown 
in Figs. \ref{fig:path}, \ref{fig:xyz-path} and \ref{fig:distancetocenter} 
and a video of the experiment which can be found at \url{https://youtu.be/E4vCSJw97FQ}.

\begin{figure}[!ht]
\centering
      \includegraphics{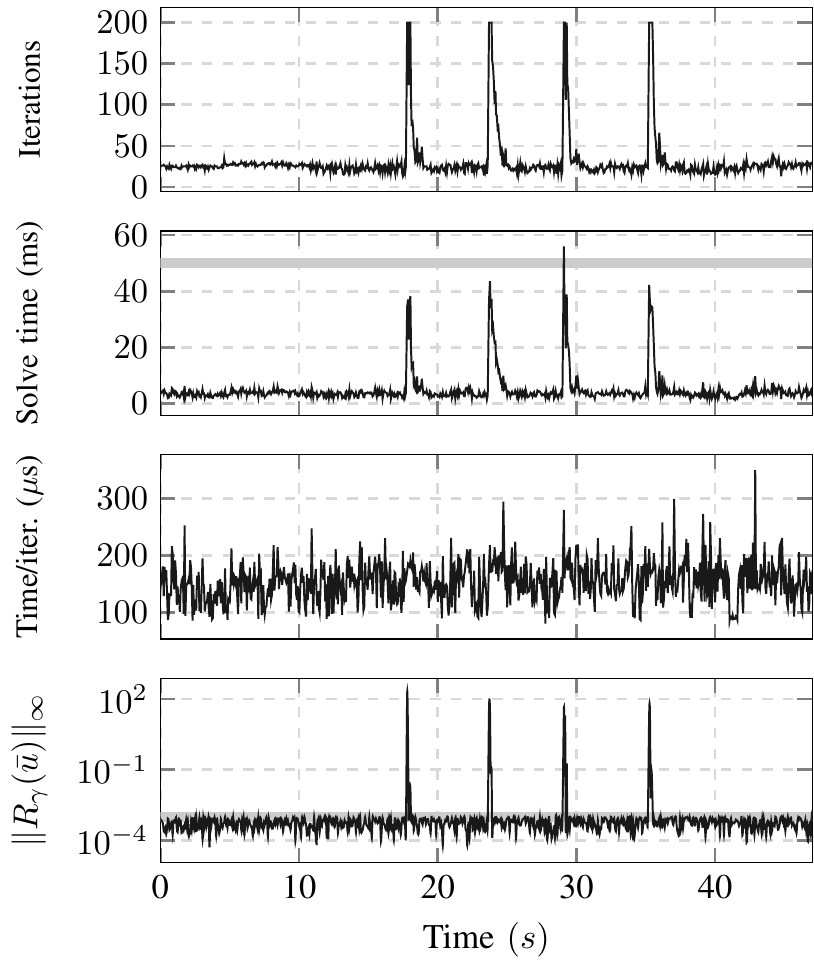}
      \caption{Solver diagnostics: (Top) Number of iterations required for
      convergence. Observe that at reference changes, the initial guess is rather
      inaccurate and the solver requires more iterations, (Middle-top) time required by PANOC to find an optimal sequence
      of control actions, (Middle-bottom) average time taken per internal iteration, and (Bottom) infinity norm of the
      fixed-point residual, \(\|R_\gamma(\bar{u})\|_{\infty}\), which serves as an indicator
      of the solution quality.}
      \label{fig:solverstats}
\end{figure}

As shown in the second subfigure of Fig.~\ref{fig:solverstats}, once the
reference changes, the solver reaches the maximum number of iterations (200 iterations)
and the solution it returns is of poor quality (fourth subfigure of
Fig.~\ref{fig:solverstats}).
This happens because at each time instance, the solver is initialized with the
previously computed optimal trajectory. Upon a reference change, the initial
guess is rather far from optimal and this necessitates more iterations for
convergence.
Nonetheless, this inaccuracy is eliminated at the next time instant ---
\(\unit[0.05]{s}\) later --- where the solver is provided a good initial
estimate and converges within the prescribed tolerance, \(\epsilon = 10^{-3}\).
This way, \ac{NMPC} is executed at \(\unit[20]{Hz}\).
As shown in the third subfigure of Fig.~\ref{fig:solverstats}, at one
time instant, the solution time exceeds the maximum allowed time. This is
accommodated by delaying the dispatch of the control action by few \(\unit{ms}\)
and has no practical effect.

The infinity norm of the fixed-point residual is below \(\epsilon\) at all
time instants with the exception of four instants from the change of reference.
Lastly, the average iteration time in every time step is shown in the
third subfigure of Fig. \ref{fig:solverstats}, and ranges from
\(\unit[80]{\mu{}s}\) to \(\unit[350]{\mu{}s}\) where the variability is
because of the different number of line search iterations.

\begin{figure}[!ht]
\centering
      \includegraphics{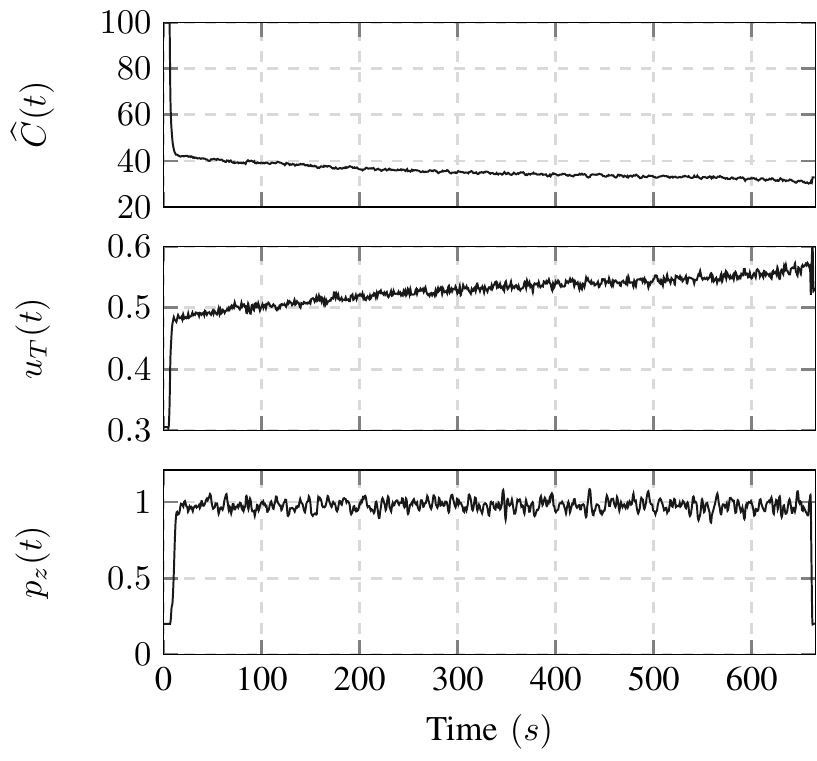}
      \caption{Control signals sent from the solver to the low-level controller during the experiment.}
      \label{fig:thrusthover}
\end{figure}

The parameters used in the dynamics of the \ac{MAV} used in the experiment are
shown in Table \ref{tab: experiment parameters}. These values were chosen empirically
(based on accurate values for other \acp{MAV}) and are not fine-tuned via experiments;
this accentuates the fact that the closed-loop and the overall obstacle
avoidance scheme is robust to errors in the determination of these parameters.

\begin{table}[!ht]
      \caption{\ac{MAV} parameters --- values in SI units}
      \label{tab: experiment parameters}
      \centering
      \begin{tabular}{c | c || c | c}
            parameter & value & parameter & value \\ \hline
            $A_x$ & 0.1 & $\tau_r$ & 0.5 \\
            $A_y$ & 0.1 & $K_r$    & 1 \\
            $A_z$ & 0.2 & $\tau_p$ & 0.5 \\
                  &  & $K_p$    & 1
      \end{tabular}
      \end{table}

The tuning parameters used by the \ac{NMPC} are
\begin{align*}
    R   {}={}& \operatorname{diag}(2,10,10)\\
    Q   {}={}& \operatorname{diag}(3I_2, 12, I_3, 3I_2),\\
    Q_f {}={}& 10 \, Q,
\end{align*}
and the prediction horizon \(T=\unit[2]{s}\).
For the \ac{EKF} for estimating the special thrust constant we have
\begin{align*}
    P_0 {}={}& 100, \\
    Q_T {}={}& 10^{-3}, \\
    R_T {}={}& 1,
\end{align*}
where $P_0$ is the initial variance for $\widehat{C}$,
$Q_T$ is the process variance in \eqref{eq:c_model}, and $R_T$ is the
measurement variance.

A separate experiment was carried out where the \ac{MAV} was given a position reference to hold for as long as
the battery could deliver power safely. This experiment was conducted to demonstrate the thrust constant estimation
described in Section~\ref{secsub:thrust estimation} and the results are presented in Fig.~\ref{fig:thrusthover}.
As the battery drains, the special thrust constant is decreasing and the control signal is adapted to
keep the \ac{MAV} hovering at a constant altitude. This experiment is part of the same video
mentioned in this section, found at \url{https://youtu.be/E4vCSJw97FQ}.

%% file: conclusions.tex
\section{Conclusions}\label{sec:conclusions}
We presented an obstacle and collision avoidance methodology coupled with 
an adaptive thrust controller that leads to increased autonomy and 
context awareness for \acp{MAV}. Obstacle avoidance is addressed with an
\ac{NMPC} controller, which is solved using PANOC --- a simple and fast 
algorithm, which involves simple algebraic operations and, unlike SQP,
does not require the solution of linear systems at each step.
Experiments were performed with the solver running onboard a \ac{MAV}
which maneuvered gently around a virtual obstacle with a smooth trajectory.
The \ac{MAV} passed the edge of the virtual obstacle
with a minimal constraint violation, as expected from the solver.

Moreover, experiments were performed to demonstrate that our thrust 
estimation method successfully compensates for the reduction of
thrust over time, making the control scheme applicable to any \ac{MAV} platform.

Future work will focus on experiments in presence of moving obstacles 
with uncertain trajectories.